\documentclass[11pt, oneside]{amsart}

\usepackage[text={5.58in,8.5in},centering,letterpaper,dvips]{geometry}
\usepackage{amssymb}
\usepackage{subcaption}
\usepackage{amsthm}
\usepackage{amsmath}
\usepackage{graphicx}
\usepackage{mathabx}
\usepackage{hyperref}
\usepackage{relsize}
\usepackage{pinlabel}
\usepackage[inline]{showlabels}
\usepackage{tikz}
\usetikzlibrary{  cd,  calc, positioning,  fit, arrows,  decorations.pathreplacing, decorations.markings,  shapes.geometric, backgrounds,bending}
\usepackage{tikzsymbols}
\theoremstyle{definition}

\theoremstyle{remark}

\begin{document}

\vspace{-.25in}
\title{An embedding into $S^2 \times S^2$}
\author{Clayton McDonald}
\address{HUN-REN Alfréd Rényi Institute of Mathematics}
\email{claytkm@gmail.com}
\urladdr{https://sites.google.com/view/claytkm/}

\vspace{-0.6 in}

\begin{abstract}
    In this note, we give a positive answer to K3 problem 4.28b, giving the first examples of homology 3-spheres that embed smoothly in $S^2 \times S^2$ but not $S^4$.
\end{abstract}
\maketitle
\vspace{-0.3 in}

A notable wrinkle to solving K3 problem 4.28b \cite{baykur2026k3} is that if a homology 3-sphere embeds in $S^2 \times S^2$, then it bounds a integer homology ball. In \cite{mcdonald2022surfacesliceshomologyspheres}, the author proves that infinitely many integer homology 3-spheres embed in a homology 4-sphere (the spin of the Poincar\'e homology sphere $P$), and do not embed in $S^4$. Therefore, all that remains to be shown is that these manifolds embed in $S^2 \times S^2$. We do this by showing that $spin(P)$ can be surgered along a curve to form $S^2 \times S^2$, and furthermore that the surgery can be done disjoint from the 3-manifold cross section. 

{\bf Part 1: Disjointness.}
The 3-manifolds in question are constructed by taking double branched covers over knots that are slices of the spin of the torus knot $T_{3,5}$. By taking double branched covers, we see a natural embedding into $spin(P)$.
The knots are even symmetric unions on $T_{3,5}$, and the double of the canonical ribbon disc for this symmetric union gives a splitting of $spin(T_{3,5})$.
%These splittings are obtained by drawing a knot on an immersed surface representing $spin(T_{3,5})$. 
%The knots split the immersed surface into two discs which by Figure 5 in \cite{mcdonald2022surfacesliceshomologyspheres} can be seen to be ribbon immersed surfaces (The self intersections of the discs are explicit in the diagram and consist of ribbon singularities).
This means that $\pi_1$ of the knot complement surjects into both discs and therefore the same is true on the level of covers.
%Therefore, $\pi_1$ of the 3-manifold surjects into both halves of the corresponding splitting of $spin(P)$. 
Given a curve represented by a homotopy class of $spin(P)$, we may then assume that it sits in the 3-manifold, and after a small push-off we may assume it is disjoint from the 3-manifold. This means we can modify $spin(P)$ away from the 3-manifold cross sections.

{\bf Part 2: The resulting manifold is $S^2 \times S^2$.}
We construct a handle decomposition for $spin(P)$ as follows: there is a balanced 2-generator presentation of $\pi_1(P)$ coming from a genus 2 Heegaard splitting of $P$, so we can represent the 4-manifold $P^{\circ} \times I$ by the corresponding 4-d 1-2 handlebody. $spin(P)$ is given by the double of this handlebody, which is obtained from the diagram of $P^{\circ} \times I$ by adding 0-framed meridians to the 2-handles, two 3-handles and a 4-handle. We now perform surgery by taking one of the dotted handles and changing it to a 0 and show that the resulting manifold is $S^2 \times S^2$. 
%First, we use the 0-framed meridians to separate the surgered handle from the rest of the diagram and cancel it with a 3-handle. 
Because this was a diagram for a homology 4-sphere, we know that the two 2-handles coming from relations in the balanced presentation link the remaining 1-handle a number of times relatively prime to one another. Therefore, we can use the Euclidean algorithm to slide them over each other and make their linking numbers with the 1-handle 0 and 1. Using the 0-framed meridians, we turn this linking number information (algebraic intersection with the 1-handle), into geometric information, and unlink the generator curves from the surgered handle. The resulting diagram is an $S^2 \times S^2$ summand, a 0-framed unknot, and a linear chain of three unknots with a dotted curve on one end and a 0 framed curve on the other. Via a 1-2 cancellation and two 2-3 cancellations, this gives the standard diagram for $S^2 \times S^2$. 
\vspace{-0.05 in}
\bibliographystyle{plain}
\bibliography{NSF}

\end{document}